# A Data Envelopment Analysis (DEA)-Based Model for Power Interruption Cost Estimation for Industrial Companies

Omid Ziaee
University of Nebraska-Lincoln
Department of Electrical and Computer Engineering
Lincoln, NE, USA
omid.ziaee@huskers.unl.edu

Bamdad Falahati
Engineering Services
Schweitzer Engineering Laboratories, Inc.
Pullman, WA, USA
bamfal@ieee.org

***Abstract*— In this paper, a new model based on Data Envelopment Analysis (DEA) and Inverse Data Envelopment Analysis (IDEA) is presented for estimating the effect of electricity on the output of industrial companies. To this end, the effect of electricity deviation, which serves as one input that can influence a manufacturing company's final product, is evaluated. Intuitively, it is known that a direct relationship exists between electricity consumption and the output of the manufacturing company.**

**However, finding a function that accurately represents this relationship is not easy. To check the applicability of the proposed method, it is tested on data from eight major vehicle manufacturing companies. In this model, labor hours, electrical energy consumption, and the value of raw materials are used as inputs, and the sales value is used as the proxy for the output of the company. These input and output data are used to find the efficiency of each company. Then, by changing the electricity consumption level, the output changes are derived. To calculate the outage cost, the deviation of the output is divided by the deviation of the electricity consumption, and the outage cost is estimated.**

## I. Introduction

The first estimates of interruption costs relied on macroeconomic indicators (i.e., gross domestic product, or wages earned) [1]-[3]. When using macroeconomic methods, a production function is estimated to model the relationship between the inputs, i.e., labor hours and electricity consumption, and the output of a specific class of industry. Then, the interruption cost can be estimated by conducting a sensitivity analysis that involves changing the electricity consumption and assessing the results on the output. In recent years, three alternative methods have replaced the traditional macroeconomic indicators method [4]-[6]. These new methods analyze customer behavior when presented over time with different reliability levels.

In the first of these methods, the market-based method, the value of uninterrupted power is considered as it relates to backup facilities (such as backup generators) and the interrupted contracts [4]. In [7], Bental and Ravid introduced a factor, namely, the ability of firms to protect against power outages by buying generators. They divided the cost of power generation into two categories; fixed cost, which is concerned with yearly capacity; and variable cost, which is concerned mainly with fuel costs. The worth of these two elements, then, is used as a proxy for interruption costs.

The second method is based on surveying selected customers. Reference [4] introduced two survey-based approaches, "direct measurement" and "direct worth". While using direct measurement, customers are asked about damages they have suffered during power interruptions, their willingness to pay (WTP) for measures to protect against power interruptions, and the minimum monetary value to compensate their damages [5]. Using the direct worth approach, outage costs are divided into two categories, the value of lost production and outage-related costs (i.e., labor costs to make up production, cost of damaged materials, etc.). Keep in mind that outage-related savings can decrease the outage costs [4]. In other words, this approach is based on estimating the production function, and it analyzes the consequences of power interruptions on production loss. As an example, in Reference [8], Tishler developed a model to estimate outage costs by using electricity and labor as exogenous variables in the production function. He assumed that the outage cost consists of three elements, namely, loss of output, reduction in productivity, and damage to materials.

The third method deals with real interruption costs. In this method, damages that occurred after a real interruption are listed to estimate outage costs. By accurately estimating customers' interruption costs, utilities can establish an efficient demand response and incentive program to enhance system reliability and resiliency [9]-[11]. In addition, in power system planning problems, using a good power interruption cost estimate leads to a better calculation of the loss load value [12].

In this paper, a new method based on Data Envelopment Analysis (DEA) and Inverse Data Envelopment Analysis (IDEA) is proposed for estimating the power interruption cost of eight major vehicle manufacturing companies. In [13]-[15], IDEA is introduced as a powerful tool for a manufacturing company to estimate its input/output level when some or all of its input/output changes. In this model, the value of raw materials, labor hours, and the amount of electrical energy consumed by a producer are selected as inputs, and the sales value as a proxy for output. By changing the amount of electrical energy and labor hours, the changes in output is calculated. Therefore, power interruption costs can be estimated using this model.

The remainder of this paper is organized as follows. In Section II, the fundamentals of DEA and IDEA are reviewed. In Section III, the relevant data and methodology of the model is introduced. Finally, the conclusion is presented in Section IV.

## II. DEA AND IDEA FUNDAMENTALS

Before explaining the mathematical formulation, in this section, the concepts of DEA and IDEA are introduced. Generally, in DEA, there are a number of producers, each of which takes some inputs to produce a set of outputs. For instance, consider the field of vehicle manufacturing. Each manufacturer uses raw materials, labor hours, and energy as inputs to produce vehicles. DEA attempts to specify which companies are most efficient and to indicate inefficiencies within the other companies.

A fundamental in this method is that if a given producer, P, can produce Y(P) units of outputs with X(P) units of inputs, then other producers should also be capable of the same production schedules if operating efficiently. This assumption should be true for a combination of producers. For example, consider a case in which two producers is combined to form a composite producer with composite inputs and composite outputs. This producer is called a "virtual producer" because it does not necessarily exist. The main goal in DEA is to find the best virtual producer for each real producer. If the virtual producer can produce more output with the same input available to the real producer, or if it can produce the same level of output with less input than the real producer, this means that the real producer is inefficient [16].

Now, consider the problem in which a producer maintains its current performance level, but some inputs increase. How much will this producer's outputs increase? This problem is considered as an inverse DEA problem because in DEA, all inputs and outputs are given, and we attempt to calculate the efficiency of each company; however, in IDEA problems, the efficiency value is given, and it is required to adjust the outputs of the producer [15].

## III. DATA AND MODEL

### A. Derivation of required data

Unlike in econometric methods, modeling by IDEA does not require time series data on inputs and outputs of the individual producers. In IDEA, if the input and output types of a set of producers are the same, i.e., similar manufacturing companies, the deviation of the output when one of the inputs changes can be estimated [13]. In this paper, all estimates are based on data gathered by Niroo Research institute (NRI), the main R&D institute of the Iranian power ministry. The dataset contains durations of outages, as well as inputs and outputs of some industrial companies in Iran.

The proposed model is applied to the vehicle manufacturing industry because this branch is one of the biggest electricity consumers, and its output is more dependent on electricity usage. As noted previously, electrical energy consumption, labor hours, and the value of raw materials as inputs, and the sales value as a proxy for output. In the interest of privacy, the real names of the vehicle manufacturing companies are not used in this paper.

### B. Model formulation

In the proposed model, the first step in estimating the power interruption cost is to reach the efficiency index for each producer. Consider n producers $P^1$, P2, …, $P^n$ , with the amount of inputs and outputs for $P^i$ being $x^i \in E^m$ and $y^i \in E^s$ , i=1,…,n, respectively. Let $x^0$ be the inputs for a producer whose efficiency and $y^0$ we want to determine as the outputs. From various DEA models, we develop the estimate based on the output-oriented CCR model (for more information, see [15]). The formulation of this model is as follows:

$$(\mathcal{P}^0)\ \max z \tag{1a}$$
$$s.t.$$
$$\textstyle\sum_{i=1}^{n} x^i \lambda_i \le x^0\ , \tag{1b}$$
$$\textstyle\sum_{i=1}^{n} y^i \lambda_i \ge z y^0, \tag{1c}$$
$$z \ge 1 \tag{1d}$$
$$\lambda_i \ge 0, \quad i = 1,2,\dots,n \tag{1e}$$

The Expression (1a) is the objective function of the problem (a) that maximizes the inverse of the efficiency for the considered company. Constraints (1b)-(1e) ensure that there is a virtual producer capable of manufacturing at least the same amount of output with less than or equal inputs. As noted in Section II, the input and output of the virtual producer are linear combinations of the other producers. Here, the inverse of z is the efficiency index of $P^0$. The optimal value of this problem must satisfy the condition $z \ge 1$ [16]. This linear optimization problem can be solved by various software, such as Excel, GAMS, MATLAB, etc.

Table I provides the relevant efficiency index of each producer. As shown, five producers are efficient, and three producers are inefficient. Now, consider an IDEA problem for $P^6$. Suppose that $P^6$ maintains its current performance, but the electricity consumption and labor hours change. How much should the sales value change?

Using the method introduced in [10], suppose the efficiency index remains unchanged, and the inputs increase from $x^0$ to $\alpha^0 = x^0 + \Delta x^0$, where $\Delta x^0 \ge 0$; it is needed to estimate the outputs $\beta = y^0 + \Delta y^0$, where $\Delta y^0 \ge 0$.

Considering the mathematical evidence in [12], if only one output exists, then the IDEA problem converts to the following model, and it is only needed to solve a single object linear programming problem:

$$(\tilde{\mathcal{P}})\ \max \beta \tag{2a}$$
$$s.t.$$
$$\textstyle\sum_{i=1}^{n} x^i \lambda_i \le \alpha^0\ , \tag{2b}$$
$$\textstyle\sum_{i=1}^{n} y^i \lambda_i \ge z^0 \beta, \tag{2c}$$

TABLE I: GATEHRED DATA FOR 8 COMPANIES AND THEIR RELATIVE EFFICIENCY

| | | $P^1$ | $P^2$ | $P^3$ | $P^4$ | $P^5$ | $P^6$ | $P^7$ | $P^8$ |
|---|---|---|---|---|---|---|---|---|---|
| Inputs | Electricity consumption (MWh) | 3174 | 14904 | 6308 | 32364 | 10866 | 161913 | 3954 | 14346 |
| | Value of raw materials ($\times 10^{10}$ Rails) | 6.55 | 52.8 | 23.0 | 294 | 65.8 | 602 | 28.7 | 6.81 |
| | Labor hours ($\times 10^6$) | 3.76 | 17.27 | 3.03 | 24.37 | 22.59 | 104.64 | 1.91 | 3.94 |
| Output | Sales value ($\times 10^{10}$ Rails) | 14.4 | 83.11 | 41.21 | 471.01 | 96.58 | 968.72 | 45.06 | 22.45 |
| Efficiency index | | 1 | 0.83914 | 1 | 1 | 0.86901 | 0.89154 | 1 | 1 |

$$\beta \geq y^0, \tag{2d}$$
$$\lambda_i \geq 0, \quad i = 1,2,\dots,n \tag{2e}$$

Where $\alpha^0$ is defined as $\alpha^0 = x^0 + \Delta x^0$, and $z^0$ is the optimal value of problem $\mathcal{P}^0$ , i.e., formulation (1a) through (1d).

Here, the change in the outputs by increasing the level of electricity consumption and labor hours is estimated. The additional amount of electricity consumption is selected by assuming the duration of the interruption and the average amount of electrical energy consumed by each producer. Considering uncertainty, a normal distribution function for uninterrupted energy consumed by a producer is supposed. Each normal distribution is introduced next by its expected value and standard deviation.

The expected value (or mean) of uninterrupted energy for each producer is calculated as follows:

$$\boldsymbol{E_{upi}} = E_{ti} + \sum_{j=1}^{m} h_{ji} D_{ji} \tag{3}$$

Where $\mathbf{E}_{upi}$ is the expected value of uninterrupted energy of producer i, $E_{ti}$ is the total amount of energy consumed by producer i, $m$ is the total number of times that producer i experiences an interruption, $h_{ji}$ is the duration of the $jth$ interruption for producer $i$, and $D_{ji}$ is the demand of producer $i$ at the start of the $jth$ interruption. Because any valid data about the demand of producers at the start of each interruption is not available, $D_{ji}$ is assumed to be the average power of each producer in a year. In other words, to calculate the $D_{ji}$, the total amount of energy consumed by each producer is divided by the total work hours of that producer. The standard deviation of this distribution is:

$$\sigma = \frac{1}{4} \times (\sum_{j=1}^{m} h_j D_j) \tag{4}$$

Taking the above assumption into consideration, the proposed model for each of the eight manufacturing companies is as follows:

$$(\tilde{\mathcal{P}}) \max SV \tag{5a}$$

$$s.t.$$
$$\textstyle\sum_{i=1}^{8} E_t^i \lambda_i \leq E^0 \quad where\ E^0 \sim N(E_{upi}, \sigma), \tag{5b}$$
$$\textstyle\sum_{i=1}^{8} L_t^i \lambda_i \leq L^0 \quad where\ L^0 = L_t^0 \times \left(\frac{E^0}{E_t^0}\right) \tag{5c}$$
$$\textstyle\sum_{i=1}^{8} R_t^i \lambda_i \leq R_t^0 \tag{5d}$$
$$\textstyle\sum_{i=1}^{8} SV_t^i \lambda_i \geq z^0 SV \tag{5f}$$
$$SV \geq SV_t^0 \tag{5g}$$
$$\lambda_i \geq 0 \quad i = 1,2,\dots,n \tag{5h}$$

Where $E_t$, $L_t$, $R_t$, and $SV_t$ are the total electricity consumption, total labor hours, total value of raw materials, and sales value for each producer, respectively; these are introduced in Table I. $R_t^0$ and $E_t^0$ represent the inputs, and $SV_t^0$ is the sales value (output) of the considered producer whose new sales value ($SV$) we want to determine by solving this optimization problem.

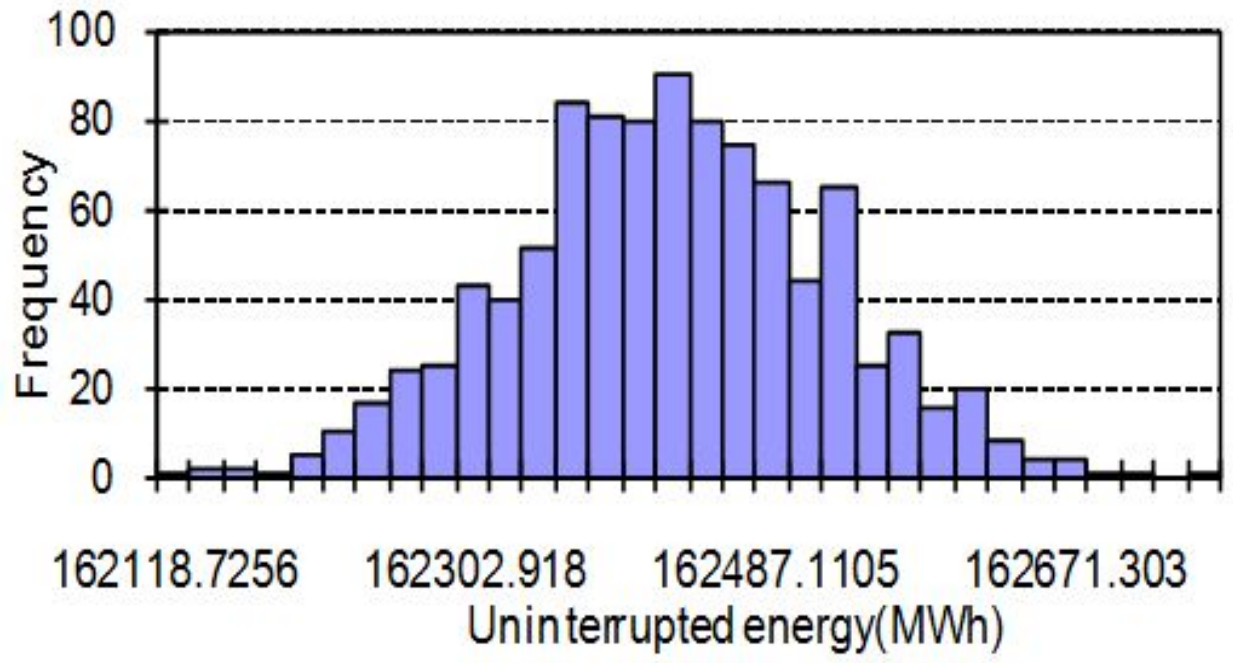

Fig. 1. Assumed distribution of uninterrupted energy for producer 6

For each producer, 1000 samples are selected randomly from the defined distribution as the new electricity consumption level ($E^0$ in problem $\tilde{\mathcal{P}}$). In addition, a linear relationship between consumed energy and labor hours is assumed. Problem $\tilde{\mathcal{P}}$ is solved for these selected numbers over 1000 iterations. The result of problem $\tilde{\mathcal{P}}$ is the optimal output level. Fig. 1 shows the selected distribution of uninterrupted energy for producer 6.

The interruption cost, per kWh, is computed by dividing the resulting increase in expected sales value by the change in electricity use, for 1000 data samples. The interruption cost calculation is formulated in (6).

TABLE II: ASSUMED DISTRIBUTION OF UNINTERRUPTED ENERGY FOR EACH PRODUCER AND THE EXPECTED VALUE OF INTERRUPTION COST

| | $P^1$ | $P^2$ | $P^3$ | $P^4$ | $P^5$ | $P^6$ | $P^7$ | $P^8$ |
|---|---|---|---|---|---|---|---|---|
| Expected value (MWh) | 3186 | 14932 | 6352 | 32450 | 11086 | 162403 | 3971 | 14455 |
| Standard deviation (MWh) | 3 | 7.5 | 11 | 10.5 | 55 | 97.5 | 4.2 | 27.2 |
| Expected value of interruption cost (Rial/KWh) | 6704.818 | 13378.31 | 9603.08 | 13576.77 | 13854.55 | 11702 | 16446.59 | 15645.68 |

$$IC = \frac{SV - SV_t^0}{E^0 - E_t^0} \quad (6)$$

In (6), IC represents the power interruption cost, SV is the result of (5), and the other parameters are the same as in the previous formulations. Equation (6) is used to estimate the interruption cost for each producer.

Table II presents the interruption costs for all producers. These results indicate that the estimated interruption costs for these eight companies varied between 6704.818 and 16446.59 (Rial/kWh). (Rial is a unit of Iranian money.)

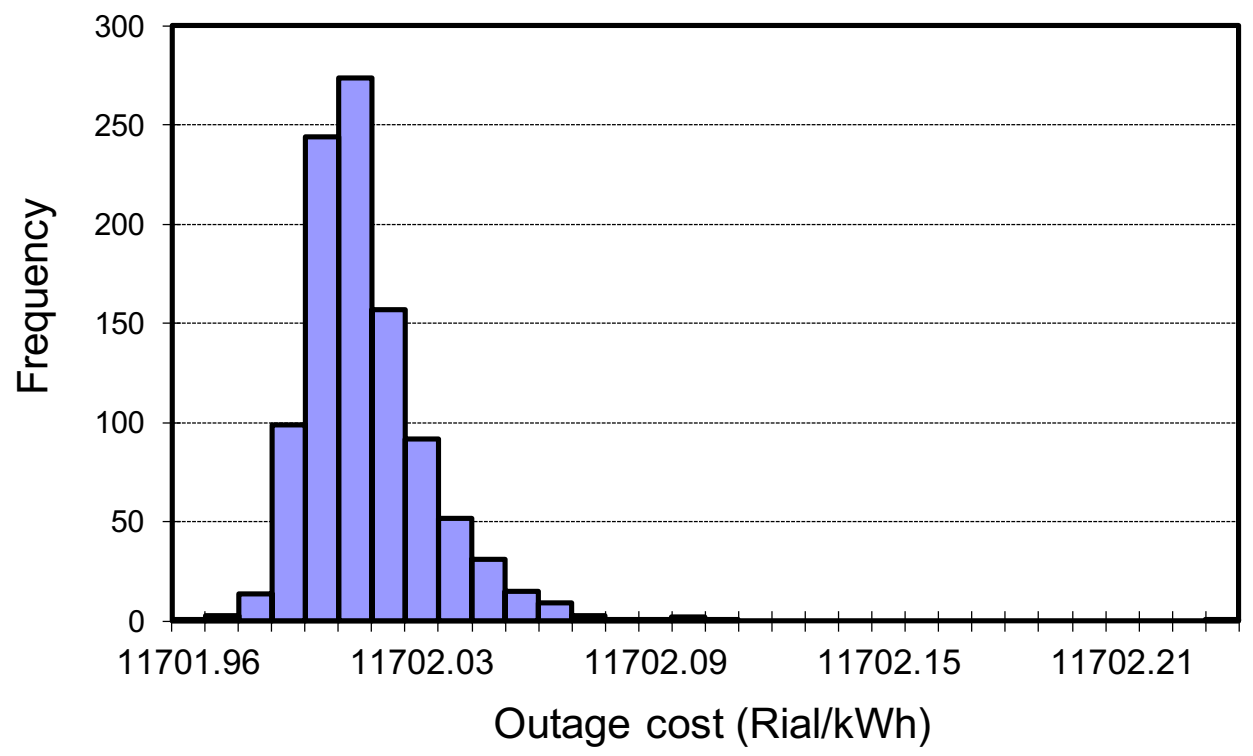

Fig. 2. Outage cost of producer 6

Fig. 2 shows the interruption cost per kWh for producer 6 for 1000 sample numbers. A lower interruption cost indicates a smaller share of electricity in production. As Table II shows, $P^1$ has the lowest and $P^7$ has the highest share of electricity in production.

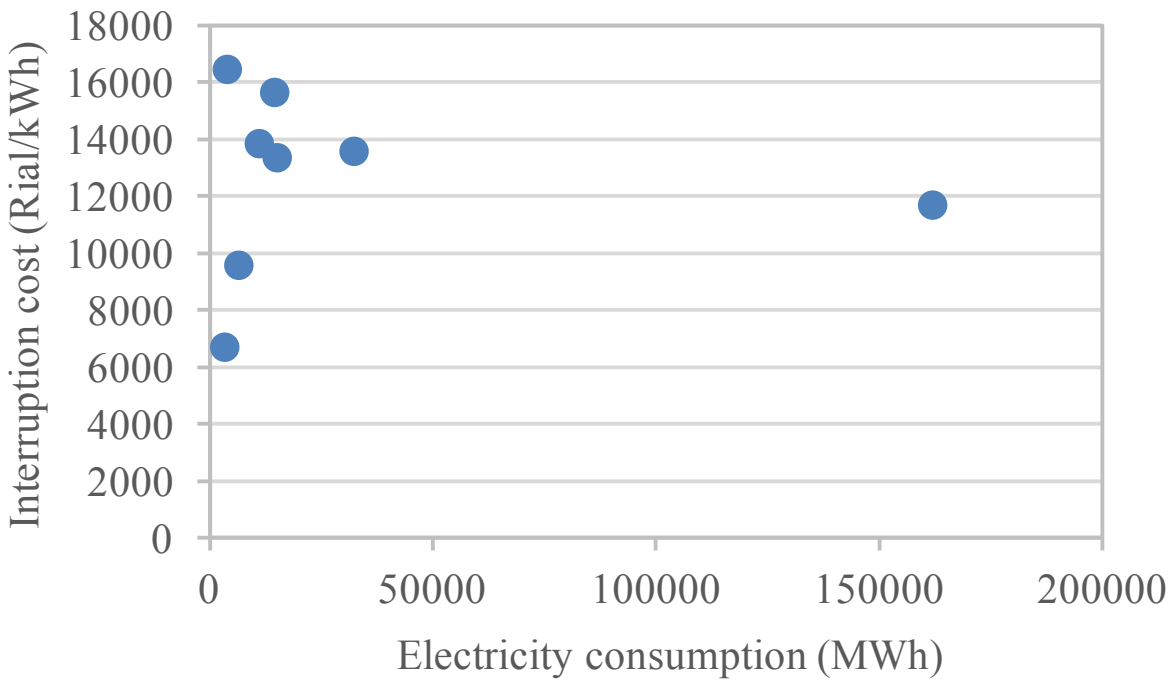

Fig. 3. Interruption cost versus electricity consumption

Fig. 3 indicates the interruption cost (vertical axis) versus the amount of electricity that each company used (horizontal axis). As can be seen in Fig. 3 no direct relationship is observed between the interruption cost and electricity consumption.

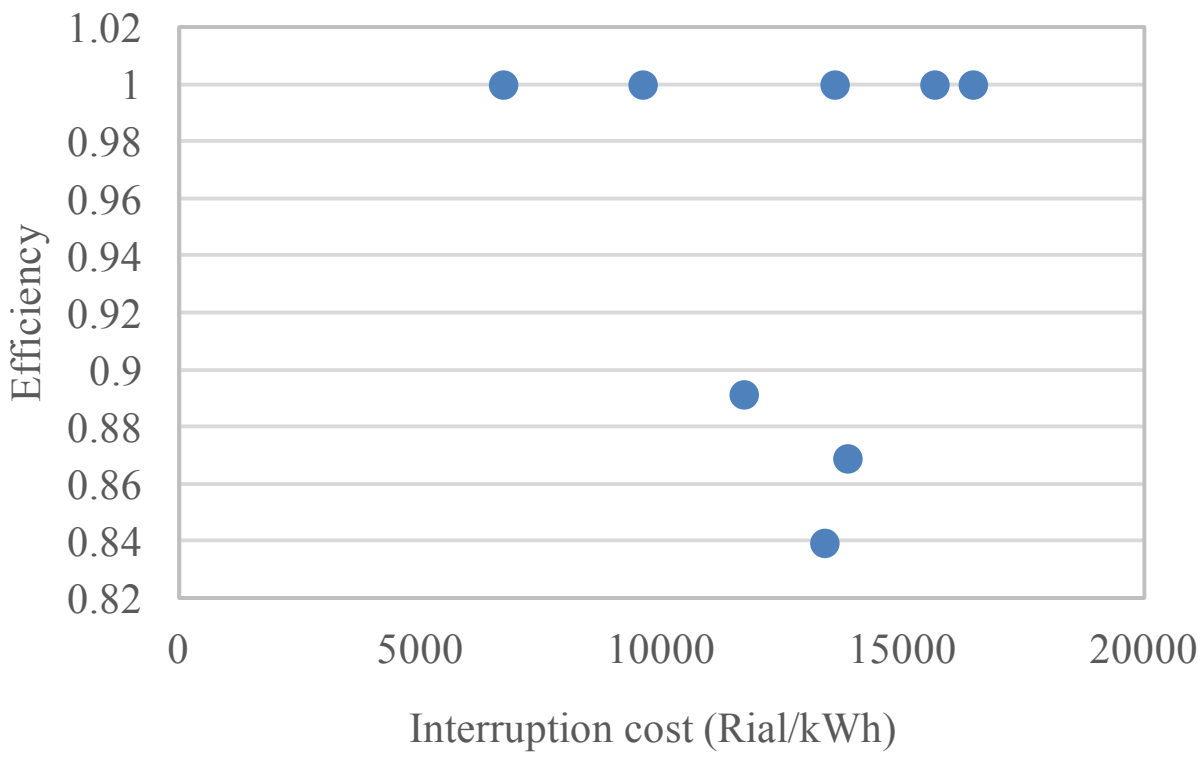

Fig. 4. Efficiency versus interruption cost for the eight producers

Fig. 4 illustrates the values of efficiency versus the interruption cost. As can be seen, no direct relationship exists between efficiency and interruption cost. This means that producers that are more sensitive to electricity outages are not always the most efficient companies. Using the other inputs efficiently also substantially influences the efficiency of the producer.

## IV. CONCLUSION

In this paper, we investigated a new approach for estimating power interruption costs based on Data Envelopment Analysis (DEA) and Inverse Data Envelopment Analysis (IDEA). Unlike econometric methods, modeling using IDEA does not require time series data on inputs and outputs of the individual producers. Using DEA, the relative efficiency of each company can be calculated.

The model was applied to eight major vehicle-manufacturing companies. The results indicate that five companies are efficient, while three companies are inefficient. To estimate the outage cost, inverse DEA is applied to these eight companies by changing the values of electricity consumption and labor hours. The results indicate that in spite of having the same inputs and outputs, the power interruption costs for these eight companies are not the same. This diversity of interruption costs shows the different share of electricity in production for these eight companies.

In short, this means that, for a company whose production is more dependent on electricity, the interruption cost is higher than for those companies whose production is less dependent on electricity consumption. In addition, the results indicate that no meaningful relationship exists between interruption cost and either electricity consumption or efficiency.

**Omid Ziaee** (S’08) earned his B.Sc. degree from Ferdowsi University, Mashhad, Iran, in 2006 and his M.Sc. degree from the University of Tehran, Tehran, Iran, in 2008. He is currently working toward his Ph.D. in electrical engineering at the University of Nebraska-Lincoln, Lincoln, NE, USA. His research interests include power systems, electricity markets, and the application of optimization and operations research to power systems. He is a lifetime member of the IEEE HKN Honor Society.

**Bamdad Falahati** (S’8, M’13) is a protection engineer with Schweitzer Engineering Laboratories (SEL). He earned his B.S. and M.S. degrees in electrical engineering from Sharif University of Technology in 1999 and 2008, respectively. He earned his Ph.D. in electrical engineering from Mississippi State University in 2013. His professional interests include power system protection, substation automation systems, power systems reliability, and distribution grid management.